\documentclass{amsart}

\usepackage[latin1]{inputenc}
\usepackage[english]{babel}
\usepackage{indentfirst}
\usepackage{amssymb}
\usepackage{amsthm}
\usepackage{xcolor}
\usepackage[all]{xy}

\newcommand{\sub}{\subseteq}

\def\epsilon{\varepsilon}

\newtheorem{theo}{Theorem}[section]
\newtheorem{cor}[theo]{Corollary}

\numberwithin{equation}{section}

\title{On weak convergence in K\"{o}the-Bochner function spaces}

\author{Jos\'{e} Rodr\'{i}guez}

\address{Dpto. de Matem\'{a}ticas\\E.T.S. de Ingenier\'{i}a Agron\'{o}mica y de Montes y Biotecnolog\'{i}a\\
Universidad de Castilla-La Mancha\\ 02071 Albacete\\ Spain} 
\email{jose.rodriguezruiz@uclm.es}



\thanks{The research was supported by grants PID2021-122126NB-C32 
(funded by MCIN/AEI/10.13039/501100011033 and ``ERDF A way of making Europe'', EU) and 
21955/PI/22 (funded by {\em Fundaci\'on S\'eneca - ACyT Regi\'{o}n de Murcia}).}

\begin{document}

\begin{abstract}
Let $E$ be an order continuous K\"{o}the function space over a non purely atomic probability measure $\mu$ 
and let $X$ be a Banach space, with topological duals $E^*$ and $X^*$, respectively. Let
$E(X)$ and $E^*(X^*)$ be the corresponding K\"{o}the-Bochner function spaces and consider $E^*(X^*)$ as a subspace of $E(X)^*$. We prove that
if $X^*$ fails the Radon-Nikod\'{y}m property, then there is a bounded, non weakly null sequence $(f_n)$ in $E(X)$ such that
$\langle \varphi,f_n\rangle \to 0$ for every $\varphi\in E^*(X^*)$; in particular, the closed unit ball of $E^*(X^*)$ is not
a James boundary for $E(X)$. This extends a result by B. Cascales and A.J. Pallar\'{e}s [Collect. Math. 45 (1994), 263--270] on the case $E=L_1(\mu)$ and 
allows us to answer a question posed recently by S. Dwivedi [Rev. Real Acad. Cienc. Exactas Fis. Nat. Ser. A-Mat. RACSAM 120 (2026), 71].
\end{abstract}

\maketitle

\section{Introduction}

Let $X$ be a Banach space and let $(\Omega,\Sigma,\mu)$ be a probability space. Let $1\leq p <\infty$ and $1<q\leq \infty$ be conjugate exponents, that is, 
$\frac{1}{p}+\frac{1}{q}=1$. Then $L_q(\mu,X^*)$ can be identified with a subspace of $L_p(\mu,X)^*$, the duality being
$$
	\langle g,f \rangle=
	\int_\Omega \langle g(\cdot),f(\cdot)\rangle \, d\mu, \qquad
	g\in L_q(\mu,X^*), \ f\in L_p(\mu,X);
$$
see, e.g., \cite[pp. 97--98]{die-uhl-J}. The equality $L_p(\mu,X)^*=L_q(\mu,X^*)$ holds
if and only if $X^*$ has the Radon-Nikod\'{y}m property with respect to~$\mu$, \cite[p.~98, Theorem~1]{die-uhl-J}, and this is always the case whenever
$\mu$ is purely atomic, \cite[p.~62]{die-uhl-J}. It is not difficult to check that, in general, 
$L_q(\mu,X^*)$ is {\em norming} for $L_p(\mu,X)$, that is, the norm of each $f\in L_p(\mu,X)$ can be computed as
\begin{equation}\label{eqn:norming}
	\|f\|_{L_p(\mu,X)}=\sup_{g\in B_{L_q(\mu,X^*)}} \langle g,f\rangle.
\end{equation}
For $p=1$, Cascales and Pallar\'{e}s~\cite{cas-pal} proved that $X^*$ has the Radon-Nikod\'{y}m property if (and only if), for 
each probability measure~$\mu$, the set $B_{L_\infty(\mu,X^*)}$ is a {\em James boundary} for~$L_1(\mu,X)$, meaning that 
in~\eqref{eqn:norming} the supremum is actually a maximum. 
Recently, Dwivedi (see \cite[Problem~3.19]{dwi}) asked whether the fact that
$B_{L_q(\mu,X^*)}$ is a James boundary for $L_p(\mu,X)$, for an atomless~$\mu$, implies that $X^*$ has the Radon-Nikod\'{y}m property. 

The aim of this note is to answer affirmatively the previous question (thus extending the aforementioned result of~\cite{cas-pal}) in the more general setting
of K\"{o}the-Bochner function spaces. 
To state our main result, Theorem~\ref{theo:main} below, we need further terminology and background.

Given a K\"{o}the function space over a probability space $(\Omega,\Sigma,\mu)$ and a Banach space~$X$,
we denote by $E(X)$ the corresponding K\"{o}the-Bochner function space, which is the Banach space of all (equivalence classes of) strongly $\mu$-measurable functions
$f:\Omega \to X$ for which the function $\|f(\cdot)\|_X:\Omega \to \mathbb{R}$ belongs to~$E$, equipped with the norm
$\|f\|_{E(X)}:=\|\|f(\cdot)\|_X\|_E$. Clearly, $E(X)$ equals the Lebesgue-Bochner function space $L_p(\mu,X)$ when $E=L_p(\mu)$ for $1\leq p\leq \infty$. If $E$ is order continuous, 
then $E^*$ is also a K\"{o}the function space over $(\Omega,\Sigma,\mu)$ (see, e.g., \cite[p.~29]{lin-tza-2}) and, as in the case of Lebesgue-Bochner function spaces:
\begin{enumerate}
\item[(i)]$E^*(X^*)$ identifies with a norming subspace of $E(X)^*$, the duality being
$$
	\langle g,f \rangle=
	\int_\Omega \langle g(\cdot),f(\cdot)\rangle \, d\mu,
	\qquad
	g\in E^*(X^*), \ f\in E(X);
$$
\item[(ii)] the equality $E(X)^*=E^*(X^*)$ holds if and only if $X^*$ has the Radon-Nikod\'{y}m property with respect to~$\mu$.
\end{enumerate}

Let $Z$ be an arbitrary Banach space. Following~\cite{nyg4}, a set $H \sub B_{Z^*}$ is said to have the {\em Rainwater property} for~$Z$ if a 
sequence $(z_n)$ in~$Z$ is weakly null if (and only if) it is bounded and $\langle z^*,z_n\rangle \to 0$ for every $z^*\in H$. 
Such terminology is motivated by Rainwater's theorem~\cite{rai}
stating that the set of all extreme points of~$B_{Z^*}$ has this property. More generally, 
according to a result of Simons~\cite{sim}, every James boundary has the Rainwater property
(cf. \cite[Theorem~3.134]{fab-ultimo}).

Now, our main result reads as:

\begin{theo}\label{theo:main}
Let $E$ be an order continuous K\"{o}the function space over a non purely atomic probability space $(\Omega,\Sigma,\mu)$ 
and let $X$ be a Banach space. If $B_{E^*(X^*)}$ has the Rainwater property for $E(X)$, then $X^*$ has the Radon-Nikod\'{y}m property. 
\end{theo}

Hence, we get the following result which provides
an affirmative answer to Dwivedi's question when applied to $E=L_p(\mu)$ for $1\leq p <\infty$:

\begin{cor}\label{cor:main}
Let $E$ be an order continuous K\"{o}the function space over a non purely atomic probability space $(\Omega,\Sigma,\mu)$ 
and let $X$ be a Banach space. If $B_{E^*(X^*)}$ is a James boundary for $E(X)$, then $X^*$ has the Radon-Nikod\'{y}m property. 
\end{cor}

Section~\ref{section:pre} introduces basic terminology and contains some preliminaries needed to deal with the proof of Theorem~\ref{theo:main}, which 
is given in Section~\ref{section:proof}. From the technical point of view, we will follow
the approach of~\cite{cas-pal}, which in turn is based on some ideas of Edgar~\cite{edg6} and Ghoussoub and Saab~\cite{gho-saa} 
involving Stegall's universal non Asplund operator~\cite{ste2}.

\section{Terminology and preliminaries}\label{section:pre}

We denote by~$\omega=\{0,1,2,\dots\}$ the set of all natural numbers and
we identify each $n\in \omega$ with the set of its predecessors, that is, $n=\{k\in \omega: k<n\}$.

All our Banach spaces are real. 
Given a Banach space~$Z$, we denote its norm by $\|\cdot\|_Z$ and its closed unit ball by $B_Z$. The density character of~$Z$, denoted by~${\rm dens}(Z)$,
is the smallest cardinality of a dense subset of~$Z$. The topological dual of~$Z$
is denoted by~$Z^*$. The evaluation of $z^*\in Z^*$ at $z\in Z$ is denoted by $\langle z^*,z\rangle$. A set $H \sub B_{Z^*}$ is said to
be a {\em James boundary} for~$Z$ if for every $z\in Z$ there is $z^*\in H$ such that $\|z\|_Z=\langle z^*,z\rangle$.
By a {\em subspace} of~$Z$ we mean a closed linear subspace. An {\em operator} is a continuous linear map between Banach spaces. 
Given a probability space $(\Omega,\Sigma,\mu)$, a function $f:\Omega \to Z$ is said to be: (i)~{\em simple} if
it can be written as $f=\sum_{i=0}^p z_i\chi_{A_i}$ for some $A_0,\dots,A_p \in \Sigma$
and $z_0,\dots,z_p\in Z$, where $\chi_{A_i}$ stands for the characteristic function of~$A_i$; (ii)~{\em strongly $\mu$-measurable}
if there is a sequence $(f_n)$ of simple $Z$-valued functions defined on~$\Omega$ such that $f_n(\omega) \to f(\omega)$ 
for $\mu$-a.e. $\omega \in \Omega$. For the basics of Bochner integration theory, see \cite[Chapter~II]{die-uhl-J}.

\subsection*{K\"{o}the function spaces}
A {\em K\"{o}the function space} over a probability space $(\Omega,\Sigma,\mu)$ is a 
Banach space~$E$ consisting of (equivalence classes of) $\mu$-integrable 
real-valued functions on~$\Omega$ such that: (i)~$\chi_A$ belongs to~$E$ for every $A\in \Sigma$; (ii)~if $|f|\leq |g|$ $\mu$-a.e., $f$ is a $\mu$-measurable real-valued function on~$\Omega$ and $g\in E$, then $f\in E$ and $\|f\|_E\leq \|g\|_E$. The identity maps 
$L_\infty(\mu)\to E$ and $E\to L_1(\mu)$ 
are well-defined, one-to-one operators (see, e.g., \cite[Lemma~2.7]{oka-alt}).
We say that $E$ is {\em order continuous} if $\|f_\alpha\|_E \to 0$ for every decreasing net $(f_\alpha)$ in~$E$ with $\inf f_\alpha=0$; in this case, the set
of all simple real-valued functions on~$\Omega$ is dense in~$E$ (see, e.g., \cite[Remark~2.6]{oka-alt}). 

\subsection*{The usual measure on~$\{0,1\}^\kappa$ and Maharam's theorem}
Let $(\Omega,\Sigma,\mu)$ be a probability space. For each $A\in \Sigma$, we write $\Sigma_A:=\{B\in \Sigma: \, B \sub A\}$,
so that $\Sigma_A$ is a $\sigma$-algebra on~$A$ and the restriction of~$\mu$ to~$\Sigma_A$, denoted by~$\mu_A$, is a non
negative finite measure. Given an infinite cardinal~$\kappa$, we say that~$\mu$ is
{\em homogeneous with Maharam type~$\kappa$} if ${\rm dens}(L_1(\mu_A))=\kappa$ for every $A\in \Sigma$ with $\mu(A)>0$. The basic example
of such a measure is the usual product probability measure $\lambda_\kappa$ on~$\{0,1\}^\kappa$. 
Recall that $\lambda_\kappa$ is defined on the $\sigma$-algebra $\Sigma_\kappa$ generated by all sets of the form
$$
	\Delta^\kappa_\sigma:=\{x\in \{0,1\}^\kappa: \, x(\alpha)=\sigma(\alpha) \, \text{ for all $\alpha \in F$}\}
$$
where $F \sub \kappa$ is finite and $\sigma\in \{0,1\}^F$.

Let us define an equivalence relation $\sim$ on~$\Sigma$ by $A \sim B$ if and only if $\mu(A\triangle B)=0$.
The equivalence class of an element~$A \in \Sigma$ is denoted by~$A^{\bullet}$. The set of such equivalence classes, denoted by $\Sigma/\mathcal{N(\mu)}$, 
becomes a {\em measure algebra}
with the usual Boolean algebra operations and the functional defined by $\mu^{\bullet}(A^{\bullet}):=\mu(A)$ for all $A\in \Sigma$.  
Maharam's theorem (see, e.g., \cite[Section~3]{fre14} or \cite[\S14]{lac-J}) states that if $\mu$ is homogeneous with Maharam type~$\kappa$ for some
infinite cardinal~$\kappa$, then the measure algebras of $(\Omega,\Sigma,\mu)$ and $(\{0,1\}^\kappa,\Sigma_\kappa,\lambda_\kappa)$ are isomorphic, that is,
there is a Boolean algebra isomorphism
$\theta: \Sigma_\kappa/\mathcal{N}(\lambda_\kappa) \to \Sigma/\mathcal{N}(\mu)$
such that $\mu^{\bullet}\circ \theta=\lambda_\kappa^{\bullet}$. 

\subsection*{Stegall's universal non Asplund operator}

We will need a result due to Stegall (see \cite[Theorem~4]{ste2}) connecting
the failure of the Radon-Nikod\'{y}m property in a dual Banach space with a certain factorization of
the so-called Haar operator, see Theorem~\ref{theo:Stegall} below. More recent related results
can be found in~\cite{bro2}.

Let $\Delta:=\{0,1\}^\omega$ be the {\em Cantor set} and
let $\mathcal{T}:=\bigcup_{n\in \omega}\{0,1\}^n$ be the {\em dyadic tree}, that is, the set of all finite sequences of $0$'s and $1$'s
(the empty sequence being its root). The collection
$\{\Delta^\omega_\tau:\tau\in \mathcal{T}\}$ is a basis of clopen sets for the 
usual topology on~$\Delta$. Given $n\in \omega$, $\tau\in \{0,1\}^n$ and $i\in \{0,1\}$, we define
$\tau \smallfrown i \in \{0,1\}^{n+1}$ by declaring $(\tau \smallfrown i)(k):=\tau(k)$ for every $k<n$ and 
$(\tau \smallfrown i)(n):=i$. Observe that
\begin{itemize}
\item $\Delta^\omega_\emptyset=\Delta$ and $\Delta^\omega_\tau=\Delta^\omega_{\tau \smallfrown 0}\cup \Delta^\omega_{\tau \smallfrown 1}$ 
for every $\tau\in \mathcal{T}$;
\item for each $n\in \omega$, the collection $\{\Delta^\omega_\tau: \tau\in \{0,1\}^n\}$ is a partition of~$\Delta$
into $2^n$ subsets with $\lambda_\omega(\Delta^\omega_\tau)=2^{-n}$ for every $\tau\in \{0,1\}^n$.
\end{itemize}
For each $\tau\in \mathcal{T}$, the function
$$
	h_\tau:=\chi_{\Delta^\omega_{\tau \smallfrown 0}}-\chi_{\Delta^\omega_{\tau \smallfrown 1}}
$$ 
belongs to $C(\Delta)$, the Banach space of all
continuous real-valued functions on~$\Delta$.

Let $\{e_\tau:\tau\in \mathcal{T}\}$ be the usual basis of~$\ell_1(\mathcal{T})$, the Banach space of all absolutely convergent series of real numbers 
indexed on~$\mathcal{T}$. 
The {\em Haar operator}
$$
	H: \ell_1(\mathcal{T}) \to L_\infty(\lambda_\omega)
$$
is the unique operator satisfying $H(e_\tau)=h_\tau$ for every $\tau \in \mathcal{T}$.
Note that $H$ actually takes its values in $C(\Delta)$ 
as a subspace of~$L_\infty(\lambda_\omega)$. With this terminology, Theorem~4 in \cite{ste2} yields:

\begin{theo}[Stegall]\label{theo:Stegall}
Let $X$ be a Banach space such that $X^*$ fails the Radon-Nikod\'{y}m property. Then there exist operators $R: \ell_1(\mathcal{T}) \to X$
and $S:X \to L_\infty(\lambda_\omega)$ such that $H=S\circ R$, where $H$ is the Haar operator as defined above.
\end{theo}

\section{Proof of Theorem~\ref{theo:main}}\label{section:proof}

We divide the proof into several steps.

\subsection*{Step 1} Since $\mu$ is not purely atomic, there exist an infinite cardinal~$\kappa$ and $\widetilde{\Delta} \in \Sigma$ with $\mu(\widetilde{\Delta})>0$ such that 
the probability measure $\tilde{\mu}:=\frac{1}{\mu(\widetilde{\Delta})}\mu_{\widetilde{\Delta}}$ is 
homogeneous with Maharam type~$\kappa$ (see, e.g., \cite[p.~122, Theorem~7]{lac-J}). 
Let $$\theta: \Sigma_\kappa / \mathcal{N}(\lambda_\kappa) \to \Sigma_{\widetilde{\Delta}}/\mathcal{N}(\tilde{\mu})$$
be a measure algebra isomorphism as given by Maharam's theorem (see Section~\ref{section:pre}).
For each finite set $F \sub \kappa$ and for each $\sigma\in \{0,1\}^F$, we
choose $\widetilde{\Delta}_\sigma\in \Sigma_{\widetilde{\Delta}}$ such that 
$$
	\theta((\Delta^\kappa_\sigma)^\bullet)=\widetilde{\Delta}_\sigma^\bullet.
$$
It is not difficult to check that we can assume further that:
\begin{itemize}
\item $\widetilde{\Delta}_\emptyset=\widetilde{\Delta}$ and
$\widetilde{\Delta}_\tau=\widetilde{\Delta}_{\tau \smallfrown 0}\cup \widetilde{\Delta}_{\tau \smallfrown 1}$ for every $\tau\in \mathcal{T}$;
\item for each $n\in \omega$, the collection $\{\widetilde{\Delta}_\tau: \tau\in \{0,1\}^n\}$ is a partition of~$\widetilde{\Delta}$
into $2^n$ subsets with $\mu(\widetilde{\Delta}_\tau)=2^{-n}\mu(\widetilde{\Delta})$ for every $\tau\in \{0,1\}^n$.
\end{itemize}
For each $\tau\in \mathcal{T}$ we define a simple real-valued function on~$\Omega$ by
$$
	\widetilde{h_\tau}:=\chi_{\widetilde{\Delta}_{\tau \smallfrown 0}}-\chi_{\widetilde{\Delta}_{\tau \smallfrown 1}}.
$$

\subsection*{Step~2} Given $n\in \omega$, let $g_n: \Omega \to \ell_1(\mathcal{T})$ be the simple function defined by
$$
	g_n(\omega):=\sum_{\tau \in \{0,1\}^n} \widetilde{h_\tau}(\omega) e_\tau
	\qquad\text{for every $\omega \in \Omega$}.
$$
Note that 
$$
	\|g_n(\omega)\|_{\ell_1(\mathcal{T})}=\sum_{\tau \in \{0,1\}^n} |\widetilde{h_\tau}(\omega)|=
	\begin{cases}
	 1 & \text{if $\omega \in \widetilde{\Delta}$},\\
	 0 & \text{if $\omega \in \Omega \setminus \widetilde{\Delta}$}.
	\end{cases}
$$
Hence, for each $m\in \omega$, the simple function $f_m:\Omega \to \ell_1(\mathcal{T})$ defined by 
$$
	f_m:=\frac{1}{m+1}\sum_{n\leq m} g_n
$$
satisfies 
\begin{equation}\label{eqn:fm}
	\|f_m(\omega)\|_{\ell_1(\mathcal{T})}\leq 1
	\qquad \text{for every $\omega \in \Omega$}.
\end{equation}

\subsection*{Step~3} We claim that, for each finite set $F \sub \kappa$ and for each $\sigma\in \{0,1\}^F$, we have
\begin{equation}\label{eqn:limit-basic}
	\lim_{m\to \infty}\left\|\int_{\widetilde{\Delta}_\sigma} f_m \, d\mu\right\|_{\ell_1(\mathcal{T})}= 0. 
\end{equation}	
As usual, given $J \sub I \sub \kappa$ and $\tau\in \{0,1\}^I$, we write $\tau|_J\in \{0,1\}^J$ to denote the restriction of~$\tau$ to~$J$.
In order to prove~\eqref{eqn:limit-basic}, observe first that, for each $n\in \omega$ such that $n>\max (F \cap \omega)$, we have
\begin{equation}\label{eqn:measures}
	\mu(\widetilde{\Delta}_\sigma \cap \widetilde{\Delta}_{\tau \smallfrown 0})=\mu(\widetilde{\Delta}_\sigma \cap \widetilde{\Delta}_{\tau \smallfrown 1})
	\quad
	\text{for every $\tau\in \{0,1\}^n$}.
\end{equation}
Indeed, 
for each $i\in \{0,1\}$ we have
$$
	\Delta^\kappa_\sigma \cap \Delta^\kappa_{\tau \smallfrown i}=
	\begin{cases}
	\Delta^\kappa_{\sigma|_{F\setminus \omega}} \cap \Delta^\kappa_{\tau\smallfrown i} & \text{if $\tau|_{F\cap \omega}=\sigma|_{F\cap  \omega}$}, \\
	\emptyset & \text{otherwise},
	\end{cases}
$$
so that 
$$
	\lambda_\kappa(\Delta^\kappa_\sigma \cap \Delta^\kappa_{\tau \smallfrown i})=
	\begin{cases}
	2^{-|F\setminus \omega|-n-1} & \text{if $\tau|_{F\cap \omega}=\sigma|_{F\cap  \omega}$}, \\
	0 & \text{otherwise},
	\end{cases}
$$
which clearly implies equality~\eqref{eqn:measures}. (Here $|F\setminus \omega|$ stands for the cardinality of~$F\setminus \omega$.)

Now, choose $k\in \omega$ large enough such that $k>\max(F\cap \omega)$. 
For every $m\geq k$ we have
$$
	\int_{\widetilde{\Delta}_\sigma} f_m\, d\mu=
	\frac{1}{m+1}\sum_{n\leq m} \sum_{\tau\in \{0,1\}^n} \big(\mu(\widetilde{\Delta}_\sigma\cap \widetilde{\Delta}_{\tau \smallfrown 0})-\mu(\widetilde{\Delta}_\sigma\cap \widetilde{\Delta}_{\tau \smallfrown 1})\big) e_\tau
$$
and therefore
\begin{eqnarray*}
	\left\|\int_{\widetilde{\Delta}_\sigma} f_m\, d\mu\right\|_{\ell_1(\mathcal{T})} &=&
	\frac{1}{m+1}\sum_{n< k} \sum_{\tau\in \{0,1\}^n} \big|\mu(\widetilde{\Delta}_\sigma\cap \widetilde{\Delta}_{\tau \smallfrown 0})-\mu(\widetilde{\Delta}_\sigma\cap \widetilde{\Delta}_{\tau \smallfrown 1})\big|
	\\ &+&
	\frac{1}{m+1}\sum_{k\leq n \leq m} \sum_{\tau\in \{0,1\}^n} \big|\mu(\widetilde{\Delta}_\sigma\cap \widetilde{\Delta}_{\tau \smallfrown 0})-\mu(\widetilde{\Delta}_\sigma\cap \widetilde{\Delta}_{\tau \smallfrown 1})\big| \\
	& \stackrel{\eqref{eqn:measures}}{=} &
	\frac{1}{m+1}\sum_{n< k} \sum_{\tau\in \{0,1\}^n} \big|\mu(\widetilde{\Delta}_\sigma\cap \widetilde{\Delta}_{\tau \smallfrown 0})-\mu(\widetilde{\Delta}_\sigma\cap \widetilde{\Delta}_{\tau \smallfrown 1})\big|\\
	& \leq &
	\frac{1}{m+1}\sum_{n< k} \sum_{\tau\in \{0,1\}^n} \mu(\widetilde{\Delta}_\sigma\cap \widetilde{\Delta}_{\tau}) \\
	& = &
	\frac{k \mu(\widetilde{\Delta}_\sigma)}{m+1}.
\end{eqnarray*}
From the last inequality it follows that~\eqref{eqn:limit-basic} holds, as claimed.

\subsection*{Step~4} We claim that, for each $A \in \Sigma$, we have 
\begin{equation}\label{eqn:limit-general}
	\lim_{m\to \infty}\left\|\int_A f_m \, d\mu\right\|_{\ell_1(\mathcal{T})}= 0.
\end{equation}
Indeed, fix $\epsilon>0$. Let $B \in \Sigma_\kappa$ such that $\theta(B^\bullet)=(A\cap \widetilde{\Delta})^\bullet$.
Then there exist finite sets $F_0,\dots,F_n \sub \kappa$ and maps $\sigma_i\in \{0,1\}^{F_i}$ such that
$$
	\lambda_\kappa \left(B\triangle \bigcup_{i\leq n} \Delta^\kappa_{\sigma_i}\right)\leq \frac{\epsilon}{\mu(\widetilde{\Delta})}
$$ 
(see, e.g., \cite[254F(e)]{freMT-2}). 
Moreover, the $\Delta^\kappa_{\sigma_i}$'s can be taken pairwise disjoint, as it is not difficult to check. 
Define $A_0:=\bigcup_{i\leq n} \widetilde{\Delta}_{\sigma_i}\in \Sigma_{\widetilde{\Delta}}$, so that 
$$
	\theta\left(\left(\bigcup_{i\leq n} \Delta^\kappa_{\sigma_i}\right)^\bullet\right)=A_0^\bullet
	\qquad\text{and}\qquad 
	\mu((A\cap \widetilde{\Delta})\triangle A_0)\leq \epsilon.
$$

On the one hand, since $\mu(\widetilde{\Delta}_{\sigma_i} \cap \widetilde{\Delta}_{\sigma_j})=0$ 
(because $\Delta^\kappa_{\sigma_i}\cap \Delta^\kappa_{\sigma_j}=\emptyset$) whenever $i\neq j$, we have
$$
	\int_{A_0}f_m \, d\mu=
	\sum_{i\leq n} \int_{\widetilde{\Delta}_{\sigma_i}} f_m \, d\mu
	\qquad
	\text{for every $m\in \omega$}
$$
and so {\em Step~3} ensures the existence of $m_0 \in \omega$ such that
\begin{equation}\label{eqn:limit-A0}
	\left\|\int_{A_0}f_m \, d\mu\right\|_{\ell_1(\mathcal{T})}\leq \epsilon
	\qquad
	\text{for every $m>m_0$}.
\end{equation}

On the other hand, since $\mu((A\cap \widetilde{\Delta}) \triangle A_0)\leq \epsilon$, 
inequality~\eqref{eqn:fm} yields
\begin{equation}\label{eqn:MORE}
	\left\|\int_{(A\cap\widetilde{\Delta})\setminus A_0} f_m\, d\mu\right\|_{\ell_1(\mathcal{T})} \leq \epsilon
	\qquad\text{and}\qquad
	\left\|\int_{A_0 \setminus (A\cap\widetilde{\Delta})} f_m\, d\mu\right\|_{\ell_1(\mathcal{T})} \leq \epsilon
\end{equation}
for every $m\in \omega$. Hence, we have
\begin{eqnarray*}
	\left\|\int_A f_m\, d\mu\right\|_{\ell_1(\mathcal{T})} &=&
	\left\|\int_{A\cap\widetilde{\Delta}} f_m\, d\mu\right\|_{\ell_1(\mathcal{T})} \\
	&\leq& 
	\left\|\int_{A_0} f_m\, d\mu\right\|_{\ell_1(\mathcal{T})} \\ &+&
	\left\|\int_{(A\cap\widetilde{\Delta})\setminus A_0} f_m\, d\mu\right\|_{\ell_1(\mathcal{T})} \\ &+&\left\|\int_{A_0\setminus (A\cap\widetilde{\Delta})} f_m\, d\mu\right\|_{\ell_1(\mathcal{T})}\ \ \stackrel{\eqref{eqn:limit-A0}\, \& \, \eqref{eqn:MORE}}{\leq} \ \ 3\epsilon
\end{eqnarray*}
for every $m>m_0$. As $\epsilon>0$ is arbitrary, \eqref{eqn:limit-general} holds.

\subsection*{Step~5} Suppose that $X^*$ fails the Radon-Nikod\'{y}m property. By Stegall's Theorem~\ref{theo:Stegall}, there 
exist operators $R: \ell_1(\mathcal{T}) \to X$ and $S:X \to L_\infty(\lambda_\omega)$ such that $H=S\circ R$.
Then we have a commutative diagram
$$
	\xymatrix@R=2pc@C=2pc{
	E(\ell_1(\mathcal{T})) \ar[rr]^{\widetilde{H}} \ar[dr]_{\widetilde{R}} &      & E(L_\infty(\lambda_\omega)) \\
	 & E(X) \ar[ur]_{\widetilde{S}} &
	}
$$
where $\widetilde{H}$, $\widetilde{R}$ and $\widetilde{S}$ are the natural composition operators. 

For each $m\in \omega$ we have $f_m\in E(\ell_1(\mathcal{T}))$ (it is a simple function) and we define $\widetilde{f_m}:=\widetilde{R}(f_m)=R\circ f_m$
(which is simple, too). 
We will show that the sequence 
$(\widetilde{f_m})$ in~$E(X)$ satisfies the following properties:
\begin{enumerate}
\item[(i)] $(\widetilde{f_m})$ is bounded.
\item[(ii)] $\langle \varphi, \widetilde{f_m} \rangle \to 0$ for every $\varphi \in E^*(X^*)$.
\item[(iii)] $(\widetilde{f_m})$ is not weakly null.
\end{enumerate}

\subsection*{Step~6} Observe that the sequence $(f_m)$ is bounded in $E(\ell_1(\mathcal{T}))$, by~\eqref{eqn:fm}
and the continuity of the identity map from $L_\infty(\mu)$ to~$E$. Hence, the sequence $(\widetilde{f_m})$ is bounded in $E(X)$.

\subsection*{Step~7} Fix $\varphi \in E^*(X^*)$. We will show that $\langle \varphi,\widetilde{f_m}\rangle \to 0$. 

Suppose first that $\varphi$ is simple and write $\varphi=\sum_{i=0}^p x_i^*\chi_{A_i}$ for some $x_i^*\in X^*$ and $A_i\in \Sigma$. In this case, 
for each $m\in \omega$ we have
$$
	\langle \varphi, \widetilde{f_m} \rangle= \sum_{i=0}^p
	\int_{A_i} x_i^* \circ \widetilde{f_m} \, d\mu=\sum_{i=0}^p
	\int_{A_i} (x_i^*\circ R) \circ f_m \, d\mu=\sum_{i=0}^p
	\left\langle x_i^* \circ R,\int_{A_i} f_m\, d\mu\right\rangle.
$$
From~\eqref{eqn:limit-general} applied to each~$A_i$ it follows that $\langle \varphi,\widetilde{f_m}\rangle \to 0$.

Let us turn to the general case.
Fix $\epsilon>0$. Since $\varphi$ is Bochner $\mu$-integrable
(bear in mind that $E^*(X^*)$ is contained in $L_1(\mu,X^*)$), there is $\delta>0$ such that 
$\int_B \|\varphi(\cdot)\|_{X^*}\, d\mu\leq \epsilon$ for every $B\in \Sigma$ with $\mu(B)\leq \delta$. 
Now, the strong $\mu$-measurability of~$\varphi$ ensures the existence of a simple function $\varphi_0:\Omega \to X^*$ 
and $A\in \Sigma$ with $\mu(\Omega \setminus A)\leq \delta$ such that $\|\varphi(\omega)-\varphi_0(\omega)\|_{X^*}\leq \epsilon$
for every $\omega \in A$ and $\varphi_0(\omega)=0$ for every $\omega \in \Omega \setminus A$ (see, e.g., \cite[p.~42, Corollary~3]{die-uhl-J}). 

On the one hand, since $\varphi_0$ is simple, we can take $m_0\in \omega$ such that 
\begin{equation}\label{eqn:aprox7}
	|\langle \varphi_0,\widetilde{f_m}\rangle|\leq \epsilon
	\quad\text{for every $m>m_0$}.
\end{equation}
On the other hand, for each $m\in \omega$ we have
\begin{eqnarray}\label{eqn:aprox77}
	|\langle \varphi-\varphi_0,\widetilde{f_m}\rangle| &\leq& \int_A \big|\langle \varphi(\cdot)-\varphi_0(\cdot),\widetilde{f_m}(\cdot) \rangle\big|\, d\mu
	+\int_{\Omega\setminus A} \big|\langle \varphi(\cdot),\widetilde{f_m}(\cdot) \rangle\big|\, d\mu \nonumber \\
	& \stackrel{\eqref{eqn:fm}}{\leq} & 
	\|R\|\int_A \|\varphi(\cdot)-\varphi_0(\cdot)\|_{X^*}\, d\mu
	+\|R\|\int_{\Omega\setminus A}\|\varphi(\cdot)\|_{X^*} \, d\mu \nonumber
	\\ &\leq& 2\|R\|\epsilon.
\end{eqnarray}
Inequalities \eqref{eqn:aprox7} and~\eqref{eqn:aprox77} yield
$$
	|\langle \varphi,\widetilde{f_m}\rangle| \leq
	|\langle \varphi_0,\widetilde{f_m}\rangle|+
	|\langle \varphi-\varphi_0,\widetilde{f_m}\rangle| \leq
	(1+2\|R\|)\epsilon
	\quad\text{for every $m>m_0$}.
$$
As $\epsilon>0$ is arbitrary, this proves that $\langle \varphi,\widetilde{f_m}\rangle\to 0$.

\subsection*{Step~8} For each $\omega \in \widetilde{\Delta}$
there is a unique $p(\omega)\in \Delta$ such that $\omega \in \widetilde{\Delta}_{p(\omega)|_n}$
for all $n\in \omega$. Clearly, we have
\begin{equation}\label{eqn:deltas}
	p^{-1}(\Delta^\omega_\tau)=\widetilde{\Delta}_\tau
	\quad\text{for every $\tau\in \mathcal{T}$}.
\end{equation}
Since every open subset of~$\Delta$ is the union of countably many sets of the form $\Delta^\omega_\tau$, 
the map $p:\widetilde{\Delta} \to \Delta$ is $\Sigma_{\widetilde{\Delta}}$-to-${\rm Borel}(\Delta)$ measurable.
Observe that~\eqref{eqn:deltas} also yields
\begin{equation}\label{eqn:deltas2}
	\widetilde{h_\tau}(\omega)h_\tau(p(\omega))=\chi_{\widetilde{\Delta}_\tau}(\omega)
	\quad
	\text{for every $\tau\in \mathcal{T}$ and for every $\omega \in \widetilde{\Delta}$}.
\end{equation}

For each $x\in \Delta$, let $\delta_x\in C(\Delta)^*$ be the evaluation functional at~$x$.
Let us consider the function $\Phi:\Omega \to C(\Delta)^*$ defined by
$$
	\Phi(\omega):=
	\begin{cases}
	\delta_{p(\omega)} & \text{if $\omega \in \widetilde{\Delta}$}, \\
	0 & \text{if $\omega \in \Omega \setminus \widetilde{\Delta}$}.
	\end{cases}
$$
Note that $\Phi$ is weak$^*$-scalarly $\Sigma$-measurable, that is,
for each $h\in C(\Delta)$ the real-valued function
$\langle h,\Phi(\cdot)\rangle$ is $\Sigma$-measurable, because $p$ is $\Sigma_{\widetilde{\Delta}}$-to-${\rm Borel}(\Delta)$ measurable and
$$
	\langle h, \Phi(\omega)\rangle=
	\begin{cases}
	h(p(\omega)) & \text{if $\omega \in \widetilde{\Delta}$}, \\
	0 & \text{if $\omega \in \Omega \setminus \widetilde{\Delta}$}.
	\end{cases}
$$
In addition, $\|\Phi(\cdot)\|_{C(\Delta)^*}=\chi_{\widetilde{\Delta}} \in E^*$. Hence, $\Phi$ induces an element $F_\Phi\in E(C(\Delta))^*$ via
the formula
$$
	F_\Phi(f):=\int_\Omega \langle \Phi(\cdot),f(\cdot)\rangle \, d\mu
	\qquad
	\text{for all $f\in E(C(\Delta))$}
$$
(see, e.g., \cite[p.~164]{lin-J}). 

Fix $m\in \omega$. Then $\widetilde{S}(\widetilde{f_m})=\widetilde{H}(f_m)$ belongs to~$E(C(\Delta))$
(as a subspace of $E(L_\infty(\lambda_\omega))$), because $H$ takes its values in~$C(\Delta)$. In fact, we have
$$
	\widetilde{H}(f_m)=H\circ f_m=\frac{1}{m+1}\sum_{n\leq m} H\circ g_n=
	\frac{1}{m+1}\sum_{n\leq m} \sum_{\tau\in \{0,1\}^n} \widetilde{h_\tau}(\cdot) h_\tau
$$
and so
\begin{eqnarray*}
	\langle F_\Phi,\widetilde{H}(f_m) \rangle &= &
	\frac{1}{m+1}\sum_{n\leq m} \sum_{\tau\in \{0,1\}^n} \int_\Omega \langle \Phi(\cdot),\widetilde{h_\tau}(\cdot) h_\tau\rangle \, d\mu \\
	& = &
	\frac{1}{m+1}\sum_{n\leq m} \sum_{\tau\in \{0,1\}^n} \int_{\widetilde{\Delta}} \widetilde{h_\tau}(\cdot) h_\tau(p(\cdot)) \, d\mu \\
	& \stackrel{\eqref{eqn:deltas2}}{=} &
	\frac{1}{m+1}\sum_{n\leq m} \sum_{\tau\in \{0,1\}^n}  \mu(\widetilde{\Delta}_\tau) \ \ = \ \ \mu(\widetilde{\Delta}).
\end{eqnarray*}
Hence, the sequence $(\widetilde{H}(f_m))$ is not weakly null in $E(C(\Delta))$
and, therefore, $(\widetilde{f_m})$ is not weakly null in~$E(X)$.

\smallskip
The proof of Theorem~\ref{theo:main} is now finished.
\qed
%

\subsection*{Acknowledgements} The research was supported by grants PID2021-122126NB-C32 
(funded by MCIN/AEI/10.13039/501100011033 and ``ERDF A way of making Europe'', EU) and 
21955/PI/22 (funded by {\em Fundaci\'on S\'eneca - ACyT Regi\'{o}n de Murcia}).


\def\cprime{$'$}\def\cdprime{$''$}
  \def\polhk#1{\setbox0=\hbox{#1}{\ooalign{\hidewidth
  \lower1.5ex\hbox{`}\hidewidth\crcr\unhbox0}}} \def\cprime{$'$}
\providecommand{\bysame}{\leavevmode\hbox to3em{\hrulefill}\thinspace}
\providecommand{\MR}{\relax\ifhmode\unskip\space\fi MR }
\providecommand{\MRhref}[2]{%
  \href{http://www.ams.org/mathscinet-getitem?mr=#1}{#2}
}
\providecommand{\href}[2]{#2}

\bibliographystyle{amsplain}

\end{document}